\documentclass[12pt]{amsart}

\usepackage{psfrag}
\usepackage{color}

\usepackage{graphicx,graphics}
\usepackage{fullpage,amssymb,amsfonts,amsmath,amstext,amsthm,amscd,enumerate,verbatim,tikz,booktabs,subcaption,orcidlink}
\usepackage[T1]{fontenc}
\usetikzlibrary{matrix,arrows}

\newtheorem{theorem}{Theorem}[section]

\newtheorem{lemma}[theorem]{Lemma}
\newtheorem{proposition}[theorem]{Proposition}
\newtheorem{corollary}[theorem]{Corollary}

\theoremstyle{definition}
\newtheorem{remark}[theorem]{Remark}
\theoremstyle{plain}

\newcommand{\Z}{\mathbb{Z}}

\newcommand{\D}{\mathbb{D}}
\newcommand{\R}{\mathbb{R}}
\newcommand{\N}{\mathbb{N}}

\newcommand{\diam}{\operatorname{diam}}

\newcommand{\dist}{\operatorname{dist}}

\newcommand{\Sph}{\mathbb{S}}
\newcommand{\F}{\mathcal{F}}
\newcommand{\J}{\mathcal{J}}

\newcommand{\Cs}{\mathcal{C}}

\numberwithin{equation}{section}

\begin{document}
\date{\today}

\title[Cantor's onion]{Cantor's onion: a uniformly quasiregular map whose Julia set is a Cantor set of spheres}
\author{Alastair N. Fletcher\,\orcidlink{0000-0003-1942-6928}}
\address{School of Mathematical and Statistical Sciences, Northern Illinois University, DeKalb, IL 60115-2888, USA}
\email{afletcher@niu.edu}

\author{Gaven J. Martin}
\address{Institute for Advanced Study, Massey University, Auckland, New Zealand}
\email{g.j.martin@massey.ac.nz}

\begin{abstract}
For integers $d,m\ge2$ with $1/d+1/m<1$, we construct a uniformly quasiregular map of $\Sph^3$ of degree $d^2+m^2$ whose Julia set is a Cantor set of concentric round spheres. To our knowledge, this is the first example of a uniformly quasiregular map in dimension three with infinitely many Fatou components. The map is a three-dimensional analogue of McMullen's maps $z^d+\lambda/z^m$, built from two Zorich-type power maps glued across a single shell by a slit construction. Its Julia set has Hausdorff dimension $2+s$, where $d^{-s}+m^{-s}=1$, and the restriction of the map to its Julia set is a skew product over the shift on two symbols and has topological entropy $\log(d^2+m^2)$. We also prove that uniformly quasiregular maps that expand uniformly near their Julia set, possibly with branching there, have no wandering domains and only finitely many periodic Fatou components, all of them attracting basins. The map constructed here belongs to this class.
\end{abstract}

\maketitle

\section{Introduction}

Let $f\colon\Sph^n\to\Sph^n$, $n\ge2$, be uniformly quasiregular (uqr): non-injective, with $K(f^k)\le K$ for some $K$ and all $k\ge1$. Here $K(g)=\max\{K_O(g),K_I(g)\}$ denotes the maximal dilatation. The Fatou set $\F(f)$ is the set where the iterates form a normal family, and the Julia set is $\J(f)=\Sph^n\setminus\F(f)$. The known examples in dimension $n\ge3$ have few Fatou components:
\begin{itemize}
\item Lattès-type maps have $\J(f)=\Sph^n$;
\item Mayer's power maps \cite{May97} have two superattracting basins separated by a round sphere, and more generally the power-type maps of \cite{FM20} have two superattracting basins separated by a quasisphere;
\item conformal trap constructions \cite{IM96,Pel99} have a single Fatou component with Cantor Julia set;
\item there are uqr maps whose Julia sets are wild or higher-genus Cantor sets \cite{FW,FS,FSV}. A Cantor set does not separate $\Sph^n$, so these too have connected Fatou set.
\end{itemize}
The first author constructed in \cite{Fl19} a uqr semigroup whose Julia set is a Cantor set of spheres. To our knowledge, the map constructed below is the first single uqr map with this property, and the first uqr map of $\Sph^3$ whose Fatou set has infinitely many components.

Round spheres have also appeared as Julia components in the non-uniform setting. In \cite[Theorem~1]{BFN}, Burkart, Fletcher and Nicks constructed a quasiregular map of $\R^3$ of transcendental type whose Julia set has components that are round spheres centred at $0$. These separate bounded hollow quasi-Fatou components, which are wandering domains escaping to $\infty$, in analogy with Baker's multiply connected wandering domains. That map is not uqr, and its dynamics is quite different from that considered here.

We follow the model of McMullen's maps $z\mapsto z^d+\lambda/z^m$ \cite{McM88}. For $1/d+1/m<1$ and small $|\lambda|$, their Julia set is homeomorphic to a Cantor set times a circle, and every bounded Fatou component is eventually mapped into the immediate basin of $\infty$. See \cite{DLU,QWY} for more on this family.

\begin{theorem}\label{thm:main}
Let $d,m\ge2$ be integers with $1/d+1/m<1$. There is a uqr map $f\colon\Sph^3\to\Sph^3$ of degree $d^2+m^2$ and a Cantor set $\Cs\subset\R$ such that
\[
\J(f)=\bigcup_{t\in \Cs}\{x:|x|=e^t\}.
\]
In particular $\dim_H\J(f)=2+s$, where $s\in(0,1)$ is the solution of $d^{-s}+m^{-s}=1$. Writing $\alpha=\min\Cs$ and $\beta=\max\Cs$, the Fatou components are exactly the following, and there are infinitely many of them:
\begin{itemize}
\item the immediate basin of $\infty$, $\{|x|>e^\beta\}\cup\{\infty\}$, which is invariant: $\infty$ is a fixed point of local index $d^2$, and $f^k\to\infty$ locally uniformly on it;
\item the ball $T=\{|x|<e^\alpha\}$, with $f(T)$ equal to the immediate basin of $\infty$;
\item for every bounded complementary interval $(a,b)$ of $\Cs$, the shell $\{e^a<|x|<e^b\}$.
\end{itemize}
Every component except the immediate basin of $\infty$ is strictly preperiodic.
\end{theorem}

The smallest degrees are $13$, for $(d,m)=(2,3)$ or $(3,2)$, which give different maps, and $18$, for $d=m=3$.

The dynamics on the Julia set is as one would expect. The restriction $f|_{\J(f)}$ is a skew product over the one-sided shift on two symbols, whose fibre maps are Lattès-type maps of the spheres. Consequently repelling periodic points are dense in $\J(f)$, and $f|_{\J(f)}$ has topological entropy $\log(d^2+m^2)=\log\deg f$, with an explicit mixing measure of maximal entropy (Proposition~\ref{prop:skew}).

The semigroup example of \cite{Fl19} is the uqr analogue of the rational semigroup generated by $z^2$ and $2z^2$. Its generators are power-type maps that fix $\infty$ and have no poles, and the Cantor structure of the Julia set comes from composing different generators.

A single uqr map of \emph{polynomial type}, that is, with $f^{-1}(\infty)=\{\infty\}$, cannot have such a Julia set. Suppose it did. Then $\infty\notin\J(f)$, so $\infty$ lies in a Fatou component $A_0$, and $f(A_0)=A_0$. Every component of $f^{-1}(A_0)$ is a Fatou component mapped onto $A_0$ \cite[Ch.~21]{IM2}, so it contains a preimage of $\infty$, which can only be $\infty$ itself. Hence $f^{-1}(A_0)=A_0$, so $A_0$ is completely invariant. By the blowing-up property of the Julia set \cite[Ch.~21]{IM2} (for $x\in\J(f)$ and any neighbourhood $N$ of $x$, the union $\bigcup_kf^k(N)$ omits at most finitely many points), every neighbourhood of a point of $\J(f)$ meets the backward orbit of $A_0$, which is $A_0$ itself. So $\J(f)\subset\partial A_0$. But $A_0$ is the exterior of the outermost sphere, so its boundary is a single sphere.

So a single map must combine $x^d$-type behaviour near $\infty$ with $\lambda/x^m$-type behaviour near a pole at $0$, as McMullen's maps do in the plane. The slit shell of Section~\ref{sec:shell} is exactly what glues these two behaviours into one branched cover.

The construction has three ingredients.
\begin{enumerate}
\item Two explicit maps semiconjugate by the Zorich map $h$ to similarities (Section~\ref{sec:zorich}). One is a rescaled Mayer power map near $\infty$. The other is an inverted power map near $0$.
\item A one-dimensional expanding system on heights $t=\log|x|$, whose invariant Cantor set gives $\Cs$ (Section~\ref{sec:1d}).
\item A quasiregular branched cover $\Phi$ of the middle shell onto the ball $\{|x|\le e^\alpha\}$, which matches the two explicit maps on the boundary (Section~\ref{sec:shell}).
\end{enumerate}
Uniform quasiregularity (Section~\ref{sec:uqr}) holds because every orbit meets the middle shell at most once, and away from it the dynamics is $h\circ(\text{similarity})\circ h^{-1}$.

Finally, Section~\ref{sec:nwd} considers uqr maps that expand uniformly near $\J(f)$, allowing branching on $\J(f)$. They have no wandering domains (Theorem~\ref{thm:nwd1}), and only finitely many periodic Fatou components, each an attracting or superattracting basin (Theorem~\ref{thm:periodic}). The map of Theorem~\ref{thm:main} shows that this class contains maps with infinitely many Fatou components.

\section{The Zorich map and semiconjugacies}\label{sec:zorich}

Write points of $\R^3$ as $(y,t)$ with $y\in\R^2$ and $t\in\R$. Let $\Gamma$ be the group of isometries of $\R^3$ generated by
\[
(y,t)\mapsto(-y,t),\qquad (y,t)\mapsto(y+2e_1,t),\qquad (y,t)\mapsto(y+2e_2,t).
\]
Let $h\colon\R^3\to\R^3\setminus\{0\}$ be the Zorich map of \cite[\S6.5.4]{IM2}, \cite{Ri,May97,BFN}, see \cite{FP21} for a detailed verification that it is quasiregular. We use the following standard properties.
\begin{enumerate}
\item[(Z1)] $h(y,t)=e^t\psi(y)$, where $\psi\colon\R^2\to\Sph^2$ is a branched cover, bi-Lipschitz on each unit square of the integer grid, and invariant under the planar part of $\Gamma$. The map $\psi$ induces a homeomorphism from the pillow $\R^2/\Gamma$ onto $\Sph^2$. It is equivariant for the reflections in the grid lines: $\psi(\bar y)=R\,\psi(y)$ for a reflection $R$ of $\Sph^2$, where $\bar y=(y_1,-y_2)$. The reflection $y\mapsto(-y_1,y_2)$ induces the same $R$, since it differs from $y\mapsto\bar y$ by $-\mathrm{id}$, which lies in $\Gamma$.
\item[(Z2)] $h$ is $K_h$-quasiregular, and its fibres are exactly the $\Gamma$-orbits. In particular $h$ maps each plane $\{t=\mathrm{const}\}$ onto the round sphere of radius $e^t$.
\item[(Z3)] The branch set of $h$ is $\{y\in\mathbb Z^2\}\times\R$, and $h$ maps it onto four rays from $0$, in the directions of the cone points $\psi(\mathbb Z^2)$.
\end{enumerate}
On each slab $\{a\le t\le b\}$ there is $L\ge1$ with $|v|/L\le|Dh(u)v|\le L|v|$ for almost every $u$ and all $v$, since $|h|=e^t$ and $\psi$ is bi-Lipschitz on each square of the grid. So $h$ is locally bi-Lipschitz off its branch set, and it is $L$-BLD on each slab \cite{MV}, \cite[I.4]{Ri}.

\begin{lemma}\label{lem:semi}
Let $F$ be an orientation-preserving similarity of $\R^3$ with $F\Gamma F^{-1}\subset\Gamma$. Then there is a unique continuous map $g\colon\R^3\setminus\{0\}\to\R^3\setminus\{0\}$ with $g\circ h=h\circ F$ on $\R^3$. It is a quasiregular branched cover with $K(g)\le K_h^2$. Moreover, if $F$ is a composition of such similarities, the corresponding map is the composition of the corresponding $g$'s.
\end{lemma}

\begin{proof}
This is essentially \cite[Theorem~21.4.1]{IM2}. Off the four rays $h(B_h)$, the map $g$ is locally of the form $h\circ F\circ h^{-1}$, for local inverses of $h$, and this gives the bound on $K(g)$ there. The rays are removable. Almost every line parallel to a coordinate axis misses them, so $g$ is ACL. Moreover $J_g$ is locally integrable, since $g$ is quasiregular off a closed null set and has locally bounded multiplicity, so $g\in W^{1,3}_{\rm loc}$ because $|Dg|^3\le K_h^2J_g$ almost everywhere. The last assertion follows from uniqueness.
\end{proof}

For the similarities used below, $F$ has ratio $d$ or $m$, and $g$ has degree $d^2$ or $m^2$ (Remark~\ref{rem:polar}).

\begin{remark}\label{rem:polar}
In polar form $x=e^t\psi(y)$, the similarity $F(y,t)=(dy,\,dt+c_1)$ gives
\[
g(x)=e^{c_1}\,e^{dt}\,\psi(dy),
\]
a rescaling of Mayer's power map of degree $d^2$. The similarity $F(y,t)=(m\bar y,\,-mt+c_2)$ gives
\[
g(x)=e^{c_2}\,e^{-mt}\,R\,\psi(my),
\]
which maps a neighbourhood of $0$ onto a neighbourhood of $\infty$. It is the analogue of $\lambda/x^m$, of degree $m^2$. In both cases $g$ maps spheres centred at $0$ onto spheres centred at $0$.
\end{remark}

\section{The one-dimensional model}\label{sec:1d}

\begin{figure}[ht]
\centering
\begin{tikzpicture}[scale=5]
\fill[gray!20] (1/3,0) rectangle (1/2,1);
\draw[thin] (0,0) rectangle (1,1);
\draw[very thick] (1/2,0) -- (1,1);
\draw[very thick] (0,1) -- (1/3,0);
\draw[dashed] (0,0) -- (1,1);
\node[below] at (0,0) {$\alpha$};
\node[below] at (1/3,0) {$\alpha_2$};
\node[below] at (1/2,0) {$\beta_1$};
\node[below] at (1,0) {$\beta$};
\node[left] at (0,1) {$\beta$};
\node at (5/12,0.5) {$\Delta$};
\node[right] at (0.8,0.55) {$\phi_1$};
\node[right] at (0.17,0.55) {$\phi_2$};
\draw[|-|,thick] (0,-0.17) -- node[below]{$I_2$} (1/3,-0.17);
\draw[|-|,thick] (1/2,-0.17) -- node[below]{$I_1$} (1,-0.17);
\end{tikzpicture}
\caption{The height dynamics for $(d,m)=(2,3)$: $\phi_2$ maps $I_2$ onto $[\alpha,\beta]$ decreasingly and $\phi_1$ maps $I_1$ onto $[\alpha,\beta]$ increasingly. The dashed line is the diagonal: $\beta$ is fixed by $\phi_1$, and $\phi_2(\alpha)=\beta$. Heights in the gap $\Delta$ correspond to the middle shell $M$ of Section~\ref{sec:shell}.}
\label{fig:heights}
\end{figure}
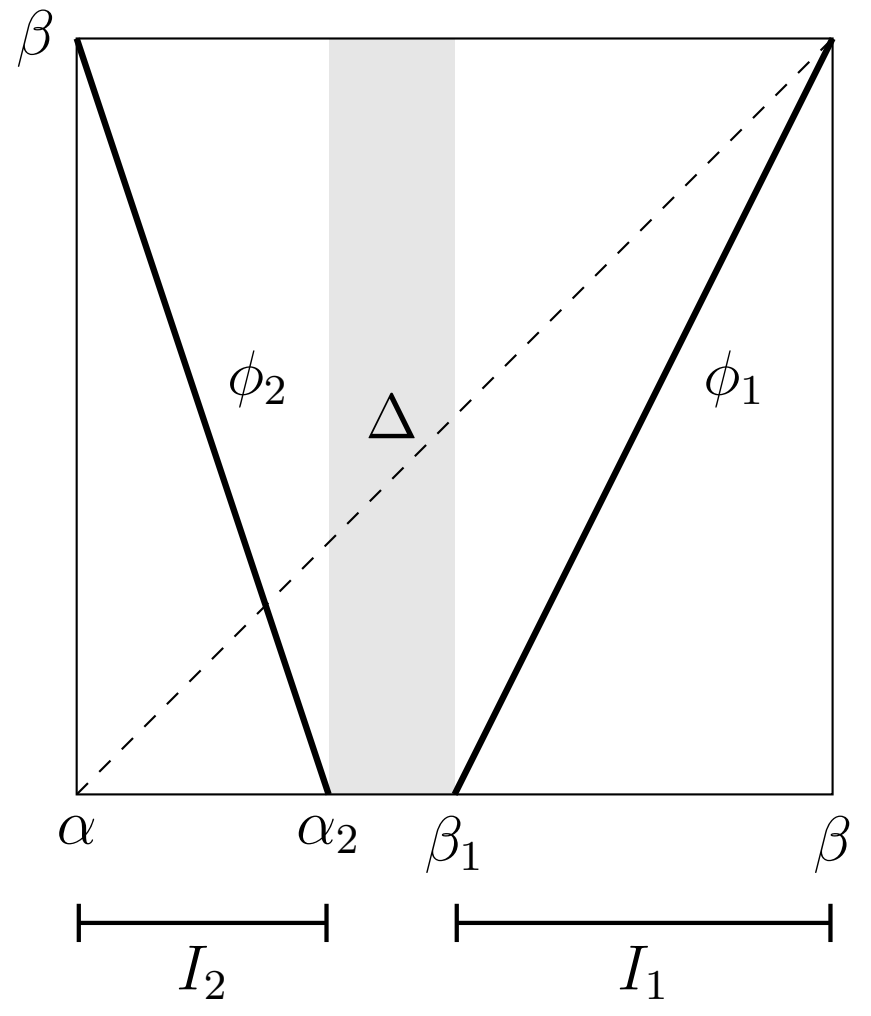

Fix $\alpha<\beta$ and put
\[
\beta_1=\beta-\frac{\beta-\alpha}{d},\qquad \alpha_2=\alpha+\frac{\beta-\alpha}{m}.
\]
Since $1/d+1/m<1$, we have $\alpha<\alpha_2<\beta_1<\beta$. Define affine maps
\[
\phi_1(t)=dt-(d-1)\beta,\qquad \phi_2(t)=-mt+\beta+m\alpha.
\]
Then
\[
\phi_1(\beta)=\beta,\quad \phi_1(\beta_1)=\alpha,\quad \phi_2(\alpha)=\beta,\quad \phi_2(\alpha_2)=\alpha.
\]
So $\phi_1$ maps $I_1=[\beta_1,\beta]$ onto $[\alpha,\beta]$ increasingly, and $\phi_2$ maps $I_2=[\alpha,\alpha_2]$ onto $[\alpha,\beta]$ decreasingly; see Figure~\ref{fig:heights}. Let $\Delta=(\alpha_2,\beta_1)$ be the gap between $I_2$ and $I_1$, and let
\[
\Cs=\{t\in I_1\cup I_2:\ \text{the orbit of $t$ under $\phi_1$ on $I_1$ and $\phi_2$ on $I_2$ stays in $I_1\cup I_2$}\}.
\]

\begin{lemma}\label{lem:1d}
$\Cs$ is a Cantor set containing $\alpha$ and $\beta$. Every $t\in[\alpha,\beta]\setminus \Cs$ has an orbit that enters $\Delta$. The bounded complementary intervals of $\Cs$ are exactly the images of $\Delta$ under compositions of the inverse branches $\phi_1^{-1}$ and $\phi_2^{-1}$; equivalently, they are the open intervals mapped onto $\Delta$ by some finite composition of the forward branches. There are infinitely many of them.
\end{lemma}

\begin{proof}
$\Cs$ is the attractor of the iterated function system $\{\phi_1^{-1},\phi_2^{-1}\}$ on $[\alpha,\beta]$. The two maps are affine contractions with ratios $1/d$ and $1/m$, and their images $I_1$ and $I_2$ are disjoint. The remaining statements are standard for such systems \cite{Hu}, \cite[Ch.~9]{Fa}. Also $\beta$ is fixed and $\phi_2(\alpha)=\beta$, so $\alpha,\beta\in \Cs$.
\end{proof}

\section{The two outer maps and the middle shell}\label{sec:shell}

Let
\[
F_1(y,t)=\big(dy,\ \phi_1(t)\big),\qquad F_2(y,t)=\big(m\bar y,\ \phi_2(t)\big).
\]
Both are orientation-preserving similarities: the linear part of $F_2$ is $\operatorname{diag}(m,-m,-m)$, with determinant $m^3>0$. Both satisfy $F_i\Gamma F_i^{-1}\subset\Gamma$. By Lemma~\ref{lem:semi} they define quasiregular maps $P$ and $Q$ of $\R^3\setminus\{0\}$ with $P\circ h=h\circ F_1$ and $Q\circ h=h\circ F_2$. Since $|P(x)|=e^{-(d-1)\beta}|x|^d$ and $|Q(x)|=e^{\beta+m\alpha}|x|^{-m}$, setting $P(0)=0$, $P(\infty)=\infty$, $Q(0)=\infty$ and $Q(\infty)=0$ gives continuous maps of $\Sph^3$, which are quasiregular because isolated points are removable.

Write $S_t=\{|x|=e^t\}$ and define the regions
\[
W_\infty=\{|x|>e^\beta\}\cup\{\infty\},\qquad T=\{|x|<e^\alpha\},\qquad M=\{e^{\alpha_2}\le|x|\le e^{\beta_1}\}.
\]
By Remark~\ref{rem:polar}, $P$ maps $S_t$ onto $S_{\phi_1(t)}$ with degree $d^2$, and $Q$ maps $S_t$ onto $S_{\phi_2(t)}$ with degree $m^2$. In particular:
\begin{itemize}
\item $P$ maps the closed ball $B_{\rm out}:=\{|x|\le e^{\beta_1}\}$ onto $\overline T$ as a proper branched cover of degree $d^2$, with $P^{-1}(\partial T)\cap B_{\rm out}=S_{\beta_1}$;
\item $Q$ maps $E_{\rm in}:=\{|x|\ge e^{\alpha_2}\}\cup\{\infty\}$ onto $\overline T$ as a proper branched cover of degree $m^2$, with $Q^{-1}(\partial T)\cap E_{\rm in}=S_{\alpha_2}$.
\end{itemize}
Note that $M=B_{\rm out}\cap E_{\rm in}$. Fix $\kappa>0$ with $\alpha_2+2\kappa<\beta_1-2\kappa$, and write $N_{\rm out}=\{e^{\beta_1-\kappa}\le|x|\le e^{\beta_1}\}$ and $N_{\rm in}=\{e^{\alpha_2}\le|x|\le e^{\alpha_2+\kappa}\}$ for the outer and inner collars of $M$.

The proposition below is consistent with degree theory. With the boundary orientation of $M$, the outer sphere $S_{\beta_1}$ is oriented as the boundary of $B_{\rm out}$ and the inner sphere $S_{\alpha_2}$ as the boundary of $E_{\rm in}$. So $P$ and $Q$ restricted to $\partial M$ map onto $\partial T$ with total degree $d^2+m^2$, which is the degree of $\Phi$.

\begin{proposition}\label{prop:shell}
There is a quasiregular branched cover $\Phi\colon M\to\overline T$ of degree $d^2+m^2$, that is, each point of $T$ has $d^2+m^2$ preimages counted with local index, such that:
\begin{itemize}
\item $\Phi=P$ on $N_{\rm out}$;
\item $\Phi=Q$ on $N_{\rm in}$;
\item $\Phi^{-1}(\partial T)=\partial M$ and $\Phi(\operatorname{int}M)=T$.
\end{itemize}
\end{proposition}

The construction of $\Phi$ is a three-dimensional version of a classical way of building branched covers of surfaces: cut sheets along slits and reglue them crosswise. In the plane, take two copies of the identity map $\D\to\D$, cut each copy along the segment $\ell=[0,\tfrac12]$, and glue the upper edge of each slit to the lower edge of the other. The result is a connected two-sheeted cover of $\D$ branched over the endpoints of $\ell$, locally modelled on $z\mapsto z^2$ near each endpoint. Preimage counts are unchanged, so the degree is the sum of the degrees of the pieces. Applied to $z^2$ on a disc and $\lambda/z^3$ on the complement of a smaller disc, over a small slit in the common target disc, this produces a degree-$5$ branched cover of an annulus onto a disc, which is the topological picture of the middle annulus of the McMullen map $z^2+\lambda/z^3$. Constructions of branched covers by cutting and pasting are standard in topology \cite[\S10.B]{Ro}, \cite{Fox}, and have a long history in the theory of quasiregular maps, going back to Rickman's work on the Picard theorem \cite{Ri85}.

The proof also uses the following relative Lipschitz approximation theorem of Tukia and Väisälä, which extends Sullivan's theorem \cite{Su}. Here a locally bi-Lipschitz embedding means an embedding that is bi-Lipschitz on a neighbourhood of each point, with respect to the Lipschitz structures.

\begin{theorem}[{\cite[Theorem~4.4]{TV}}]\label{thm:TV}
Let $M'$ and $N'$ be Lipschitz $n$-manifolds without boundary, $n\ne4$, let $A,B\subset M'$ be closed, and let $B'$ be a neighbourhood of $B$. Let $f\colon M'\to N'$ be an embedding that is a locally bi-Lipschitz embedding near $A$. Then there is an isotopy $(f_t)_{t\in[0,1]}$ of embeddings $M'\to N'$ such that:
\begin{itemize}
\item $f_0=f$;
\item $f_t=f$ on $M'\setminus B'$ and near $A$, for all $t$;
\item $f_1$ is a locally bi-Lipschitz embedding near $A\cup B$.
\end{itemize}
\end{theorem}

The theorem in \cite{TV} also controls the distance from $f_t$ to $f$, which we do not need.

\begin{proof}
\emph{Step 1: a slit.} By (Z3) and Lemma~\ref{lem:semi}, the branch set of $P$ is contained in $h(B_h\cup F_1^{-1}(B_h))\cup\{0,\infty\}$. Since $F_1(B_h)\subset B_h$, its branch values other than $0$ and $\infty$ lie on the four rays $h(B_h)$, and the same holds for $Q$. Choose a closed round ball $U\subset T$ whose centre lies off these rays, with radius small compared with the modulus of its centre. Then $U$ is at positive distance from the closures of the rays, which contain $0$. Choose also a flat closed $2$-disc $\Sigma\subset\operatorname{int}U$.

Then $P^{-1}(U)\cap B_{\rm out}$ consists of $d^2$ disjoint closed balls, each mapped homeomorphically onto $U$. Similarly $Q^{-1}(U)\cap E_{\rm in}$ consists of $m^2$ such balls. We may choose the centre of $U$ close to $0$. Then all of these balls lie away from $M$: those of $P$ in $\{|x|<e^{\alpha_2}\}$, and those of $Q$ in $\{|x|>e^{\beta_1}\}$. Fix one ball $U_P$ of the first kind and one ball $U_Q$ of the second, and put
\[
\sigma_P=(P|_{U_P})^{-1},\quad \sigma_Q=(Q|_{U_Q})^{-1},\quad D_P=\sigma_P(\Sigma),\quad D_Q=\sigma_Q(\Sigma).
\]
The maps $\sigma_P$ and $\sigma_Q$ are bi-Lipschitz. Indeed, $F_1$ and $F_2$ map $B_h=\Z^2\times\R$ into itself, so $P$ and $Q$ map $h(B_h)$ into $h(B_h)$. Since $U$ misses $h(B_h)$, so do $U_P$ and $U_Q$. On these compact sets $P=h\circ F_1\circ h^{-1}$ and $Q=h\circ F_2\circ h^{-1}$ for local inverses of $h$ away from its branch set, where $h$ is locally bi-Lipschitz.

\emph{Step 2: cross-gluing.} Let $\widehat B$ and $\widehat E$ be disjoint abstract copies of $B_{\rm out}$ and $E_{\rm in}$, on which we keep the maps $P$ and $Q$. Note that each of $\widehat B$ and $\widehat E$ contains a copy of $M$, and that the slits $D_P$ and $D_Q$ lie outside these copies. Give $\Sigma$ a transverse orientation, with sides $\Sigma^\pm$, and transport it to $D_P$ and $D_Q$ by $\sigma_P$ and $\sigma_Q$. Cut $\widehat B$ along $D_P$ and $\widehat E$ along $D_Q$, so that each slit acquires two faces, joined along its rim. Glue the $+$ face of $D_P$ to the $-$ face of $D_Q$ via $\sigma_Q\circ\sigma_P^{-1}$, and the $-$ face of $D_P$ to the $+$ face of $D_Q$ likewise. Let $X$ be the resulting space. Identified points have the same image under $P$ and $Q$, so $P$ and $Q$ together define a continuous map $G\colon X\to\overline T$.

We give $X$ an atlas with locally bi-Lipschitz transition maps, in which $G$ is quasiregular. Choose an affine coordinate on $U$ in which $\Sigma$ is a round disc in a coordinate plane. Near each open arc of the circle $\partial\Sigma$ of length less than that of $\partial\Sigma$, use smooth toroidal coordinates $c=(r,\theta,s)$, where $s$ is arc length along the arc, $(r,\theta)$ are polar coordinates in the normal plane, and $\Sigma$ corresponds to $\theta=\pi$. We read such a triple as the point $(r\cos\theta,r\sin\theta,s)$ of $\R^3$, so that the axis $r=0$ is not a coordinate singularity. There are three kinds of charts.
\begin{itemize}
\item \emph{Away from the seam}, the charts of $\Sph^3$ restricted to $\widehat B$ or $\widehat E$. There $G$ is $P$ or $Q$.
\item \emph{On the interior of the seam}, the chart is $G$ followed by the coordinate on $U$. A point just on the $+$ side lies in one piece and maps just above $\Sigma$, while a point just on the $-$ side lies in the other piece and maps just below $\Sigma$, so $G$ is a local homeomorphism there; this is what the crosswise rule is for. Its transition to the charts of the first kind is $P$ or $Q$ followed by an affine map, hence bi-Lipschitz by Step~1.
\item \emph{Along the rim} $\partial D_P\sim\partial D_Q$, a neighbourhood in $X$ consists of two copies of a neighbourhood of a point of $\partial\Sigma$ cut along $\Sigma$ and glued crosswise. In each normal plane this is the planar picture described before the proof, so the neighbourhood is the double cover of a neighbourhood of $\partial\Sigma$ in $U$ branched along $\partial\Sigma$. The chart is $\zeta=w^{-1}\circ c\circ G$, using the branch of $w^{-1}$ determined by the sheet, where $w(r,\theta,s)=(r,2\theta,s)$ is the winding map. It is a homeomorphism onto a neighbourhood in $\R^3$, and $G=c^{-1}\circ w\circ\zeta$. Its transitions to the charts above are branches of $w^{-1}$, composed with the bi-Lipschitz maps $c$, $P$ and $Q$. These overlaps avoid the axis $r=0$, where the branches of $w^{-1}$ are smooth, so the transitions are locally bi-Lipschitz. This is all that is needed for a Lipschitz atlas. Transitions between rim charts over overlapping arcs are translations in $s$, since the labelling of the sheets by $\widehat B$ and $\widehat E$ is global, so they are isometries.
\end{itemize}
In these charts $G$ is $P$, $Q$, an affine map, or $c^{-1}\circ w$. The winding map is quasiregular \cite[I.3.1]{Ri}, so $G$ is a quasiregular branched cover. Cutting and regluing does not change preimage counts, so its degree is $d^2+m^2$. Its only new branch set is the rim circle, of local index $2$, which lies over $\partial\Sigma\subset T$.

\emph{Step 3: topology.} A thin regular neighbourhood of $D_P$ in $\widehat B$ is a ball. Cutting along $D_P$ opens the slit: the two faces, joined along the rim, form a $2$-sphere bounding a spherical hole. So the cut $\widehat B$ is homeomorphic to $\widehat B$ with an open ball removed, and likewise for the cut $\widehat E$. Cross-gluing identifies the two hole-spheres, face to face and rim to rim. Removing a ball from each of two $3$-manifolds and gluing the resulting boundary spheres is the connected sum \cite[Ch.~3]{He}. Since $B_{\rm out}$ and $E_{\rm in}$ are closed balls in $\Sph^3$,
\[
X\cong B_{\rm out}\# E_{\rm in}\cong B^3\#B^3\cong S^2\times[0,1].
\]
The boundary spheres of $X$ are the copies of $S_{\beta_1}=\partial B_{\rm out}$ in $\widehat B$ and of $S_{\alpha_2}=\partial E_{\rm in}$ in $\widehat E$. The copy of $N_{\rm out}$ in $\widehat B$ and the copy of $N_{\rm in}$ in $\widehat E$ are untouched by the surgery, and on them $G=P$ and $G=Q$ respectively.

We now transfer $G$ to $M$. Both $M$ and $X$ are homeomorphic to $S^2\times[0,1]$. Orient $X$ so that $G$ is sense-preserving. This is possible: $P$ and $Q$ are sense-preserving on $\widehat B$ and $\widehat E$, and near interior points of the seam $G$ is a local homeomorphism, so the orientations of $\widehat B$ and $\widehat E$ agree across the seam. Since the rim has codimension two, they define an orientation of $X$. Call $\{e^{\beta_1-2\kappa}\le|x|\le e^{\beta_1}\}$ and $\{e^{\alpha_2}\le|x|\le e^{\alpha_2+2\kappa}\}$ the $2\kappa$-collars of $M$. Their copies in $\widehat B$ and $\widehat E$ are untouched by the surgery. Let $M_0\subset M$ and $X_0\subset X$ be the closures of the complements of the $2\kappa$-collars and of their copies. Both are again homeomorphic to $S^2\times[0,1]$. Choose an orientation-preserving homeomorphism $H\colon M_0\to X_0$ sending outer to outer and inner to inner boundary sphere; starting from any homeomorphism sending outer to outer, compose if necessary with a reflection in the $S^2$ factor. The identity maps of the $2\kappa$-collars onto their copies are orientation-preserving, since $G$ is $P$ or $Q$ there and both are sense-preserving. So on each boundary sphere of $M_0$, the map $H$ differs from the identity map onto the corresponding copy by an orientation-preserving homeomorphism of $S^2$, which is isotopic to the identity \cite{Kn}. Using these isotopies in thin collars of $\partial M_0$ inside $M_0$, we modify $H$ so that it agrees with the identity map on $\partial M_0$. Extending by the identity map of the $2\kappa$-collars onto their copies gives a homeomorphism $\Psi_0\colon M\to X$ that maps the $2\kappa$-collars identically onto their copies.

The map $\Psi_0$ maps $\partial M$ onto $\partial X$, and it is bi-Lipschitz on the $2\kappa$-collars. To make it locally bi-Lipschitz everywhere, we use Theorem~\ref{thm:TV}. We apply it with $n=3$, $M'=\operatorname{int}M$, $N'=\operatorname{int}X$, $f=\Psi_0|_{\operatorname{int}M}$, $A=\operatorname{int}M\cap(N_{\rm out}\cup N_{\rm in})$ and $B=B'=\operatorname{int}M$. Since $\Psi_0$ is bi-Lipschitz on the $2\kappa$-collars, $f$ is a locally bi-Lipschitz embedding near $A$. Let $\Psi=f_1$ on $\operatorname{int}M$ and $\Psi=\Psi_0$ on $\partial M$. Then $\Psi=\Psi_0$ on $N_{\rm out}\cup N_{\rm in}$, so $\Psi$ is continuous and locally bi-Lipschitz on all of $M$. It is injective, because $f_1(\operatorname{int}M)\subset\operatorname{int}X$ and $\Psi_0(\partial M)=\partial X$. Only a version for manifolds without boundary, as in Theorem~\ref{thm:TV}, is needed, because $\Psi_0$ is already bi-Lipschitz near $\partial M$. Since $\Psi(\partial M)\subset\partial X$, the set $\Psi(\operatorname{int}M)$ is open in $\operatorname{int}X$, by invariance of domain, and closed in $\operatorname{int}X$, by compactness. As $\operatorname{int}X$ is connected, $\Psi$ is onto. Covering $X$ by finitely many charts of the atlas above, by compactness $\Psi$ is uniformly locally bi-Lipschitz in these charts. In particular it is quasiconformal, and so $\Phi$ below has finite maximal dilatation, which is all that is used later. In the proof as given here, this approximation step is the only non-explicit part, so this proof does not give an explicit bound for $K(\Phi)$, or for the constant in Proposition~\ref{prop:uqr}. Remark~\ref{rem:highdim} outlines how $\Psi$ can instead be written down explicitly.

Set $\Phi=G\circ\Psi$. This is a quasiregular branched cover of degree $d^2+m^2$. It equals $P$ on $N_{\rm out}$ and $Q$ on $N_{\rm in}$. Also $\Phi^{-1}(\partial T)=\partial M$ and $\Phi(\operatorname{int}M)=T$, because $G^{-1}(\partial T)=\partial X$.
\end{proof}

\begin{remark}\label{rem:interp}
The boundary spheres of $M$ both cover the single sphere $\partial T$, from the same side, with partial degrees $d^2$ and $m^2$. This is a ``fold'' rather than an interpolation between distinct target spheres. The extension theorems of Berstein and Edmonds \cite{BE} and of Pankka, Rajala and Wu \cite[Theorem~3.1]{PRW} require each boundary component to cover its own target sphere with the full degree, so they do not apply directly. The quasiregular interpolations of Drasin and Sastry \cite{DS} and of Burkart, Fletcher and Nicks \cite[Theorem~2]{BFN} solve a different problem. They interpolate in a ring between power-type maps of different degrees given on the two boundary spheres, which map onto two different spheres: piecewise linear power maps in \cite{DS}, and the Mayer power maps $p_d$ and $p_{3d}$ in \cite{BFN}. Note also that only finiteness of $K(\Phi)$ matters here, since every orbit meets $M$ at most once (Proposition~\ref{prop:uqr}). The slit construction gives this directly, for every pair $(d,m)$ with $1/d+1/m<1$.
\end{remark}

\begin{remark}\label{rem:trapdoor}
The branch values of $\Phi$ lie on the four rays $h(B_h)$ and on the circle $\partial\Sigma\subset U\subset T$. Indeed, off the seam $G$ is $P$ or $Q$, and the only new branch set is the rim. So every branch value of $f$ other than $0$, $\infty$ and the points of the rays $h(B_h)$ lies in the ball $T$, which $f$ maps into the immediate basin of $\infty$. This is the analogue of the condition for the Cantor-circle regime of the McMullen maps: there, the free critical values lie in the trap door, the component of the preimage of the immediate basin of $\infty$ that contains the pole \cite{DLU}.
\end{remark}

\section{The map and its uniform quasiregularity}\label{sec:uqr}

Define $f\colon\Sph^3\to\Sph^3$ by
\[
f=\begin{cases}
P & \text{on } \{|x|\ge e^{\beta_1}\}\cup\{\infty\},\\
\Phi & \text{on } M,\\
Q & \text{on } \{|x|\le e^{\alpha_2}\}.
\end{cases}
\]
By Proposition~\ref{prop:shell}, $f=P$ on the open set $\{|x|>e^{\beta_1-\kappa}\}\cup\{\infty\}$ and $f=Q$ on the open set $\{|x|<e^{\alpha_2+\kappa}\}$. So $f$ is continuous and quasiregular. The construction depends on the choice of $\alpha<\beta$, which are arbitrary, and of $U$, $\Sigma$ and $\Psi$, so Theorem~\ref{thm:main} in fact provides a family of such maps.

\begin{lemma}\label{lem:degree}
$f$ is a branched cover of degree $d^2+m^2$. Moreover:
\begin{enumerate}
\item $f(W_\infty)\subset W_\infty$ and $f(\overline T)=\overline{W_\infty}$;
\item $f(M)=\overline T$;
\item $f$ maps $S_t$ onto $S_{\phi_1(t)}$ if $t\in I_1$, and onto $S_{\phi_2(t)}$ if $t\in I_2$.
\end{enumerate}
\end{lemma}

\begin{proof}
Count preimages, using Remark~\ref{rem:polar} and Proposition~\ref{prop:shell}.
\begin{itemize}
\item A point of $W_\infty$ has $d^2$ preimages in $W_\infty$, under $P$, and $m^2$ in $T$, under $Q$.
\item A point of $\{e^\alpha\le|x|\le e^\beta\}$ has $d^2$ preimages at heights in $I_1$ and $m^2$ at heights in $I_2$, and none in $\operatorname{int}M$, since $\Phi(\operatorname{int}M)=T$. For a point of $\partial T=S_\alpha$, these preimages lie on $\partial M$, where $f=P$ on $S_{\beta_1}$ and $f=Q$ on $S_{\alpha_2}$, so none is counted twice.
\item A point of $T$ has all $d^2+m^2$ of its preimages in $M$, under $\Phi$.
\end{itemize}
The remaining statements follow from the formulas for $\phi_1$ and $\phi_2$.
\end{proof}

\begin{proposition}\label{prop:uqr}
For all $k\ge1$, $K(f^k)\le K_h^4\,K(\Phi)$. In particular $f$ is uniformly quasiregular.
\end{proposition}

\begin{proof}
First we show that each orbit meets $M$ at most once. If $f^j(x)\in M$, then $f^{j+1}(x)\in\overline T$ and $f^{j+2}(x)\in\overline{W_\infty}$. The set $\overline{W_\infty}$ is forward invariant, and it is disjoint from $M$ because $\beta_1<\beta$. The set $\overline T$ is disjoint from $M$ because $\alpha<\alpha_2$.

Next, off $\operatorname{int}M$, $f$ is $P$ or $Q$. So along an orbit segment avoiding $\operatorname{int}M$, $f^i$ is locally $h\circ F\circ h^{-1}$, where $F$ is a composition of $F_1$ and $F_2$ and hence a similarity. By Lemma~\ref{lem:semi} its distortion is at most $K_h^2$ at a.e.\ point.

For a.e.\ $x$ and each $k$, the orbit segment $x,\dots,f^{k-1}(x)$ therefore splits into at most three blocks:
\begin{itemize}
\item a block avoiding $\operatorname{int}M$, of distortion at most $K_h^2$; it may end on $\partial M$, where $f$ is still $P$ or $Q$ on a neighbourhood;
\item at most one application of $\Phi$;
\item a block avoiding $\operatorname{int}M$ again, which after passing through $\overline T$ lies in $\overline{W_\infty}$, of distortion at most $K_h^2$.
\end{itemize}
Pointwise distortion is submultiplicative at points of differentiability outside the (null) branch sets, where the chain rule holds almost everywhere because quasiregular maps satisfy Lusin's conditions $(N)$ and $(N^{-1})$ \cite[Ch.~I]{Ri}. Submultiplicativity is used for $K_O$ and $K_I$ separately. This gives $K(f^k)\le K_h^2\cdot K(\Phi)\cdot K_h^2$.
\end{proof}

\section{Dynamics}\label{sec:dyn}

\begin{theorem}\label{thm:dyn}
We have $\J(f)=\bigcup_{t\in \Cs}S_t$, and the Fatou components are:
\begin{itemize}
\item $W_\infty$, on which $f^k\to\infty$ locally uniformly;
\item $T=\{|x|<e^\alpha\}$;
\item the shells $\{e^a<|x|<e^b\}$ for the bounded complementary intervals $(a,b)$ of $\Cs$ in $[\alpha,\beta]$.
\end{itemize}
In particular there are infinitely many Fatou components, and all except $W_\infty$ are strictly preperiodic.
\end{theorem}

\begin{proof}
On $W_\infty$, heights increase under $\phi_1$ to $+\infty$ locally uniformly, so $f^k\to\infty$ locally uniformly and $W_\infty\subset\F(f)$. Since $f(T)\subset W_\infty$, also $T\subset\F(f)$. Since $f(\operatorname{int}M)\subset T$, also $\operatorname{int}M\subset\F(f)$.

Now let $t\in[\alpha,\beta]\setminus \Cs$. By Lemma~\ref{lem:1d} and Lemma~\ref{lem:degree}(3), the sphere $S_t$, together with a shell neighbourhood of it, is mapped into $\operatorname{int}M$ after finitely many steps. Indeed, $t$ lies in a bounded complementary interval of $\Cs$, whose forward images each lie in a single $I_i$ until one of them equals $\Delta$. So $S_t\subset\F(f)$.

On the other hand, let $t\in \Cs$ and $x\in S_t$. Suppose $\{f^k\}$ were equicontinuous near $x$, with respect to the chordal metric. Put $x=h(y,t)$. For every $k$, $f^k$ maps $S_t$ onto $S_{t_k}$, where $t_k\in \Cs\subset[\alpha,\beta]$. In Zorich coordinates this map is induced by a similarity of the $y$-variable with ratio at least $2^k$.

So for any disc $D\subset S_t$ around $x$ and all large $k$, $f^k(D)$ contains the image of a full fundamental domain of $\Gamma$, that is, $f^k(D)=S_{t_k}$. These spheres have radius at least $e^\alpha$, so their diameters are bounded below. This contradicts equicontinuity. Hence $S_t\subset\J(f)$.

Each listed set in the statement of the theorem is connected, lies in $\F(f)$, and has boundary in $\J(f)$, so it is a Fatou component. These are all the components: $\Cs$ is closed with $\alpha,\beta\in\Cs$, so $[\alpha,\beta]\setminus\Cs$ is the union of the bounded complementary intervals, and the listed sets cover $\Sph^3\setminus\J(f)$. The bounded complementary intervals of $\Cs$ are infinite in number by Lemma~\ref{lem:1d}. Each shell is mapped by some $f^j$ onto $\operatorname{int}M$, since the interval is mapped onto $\Delta$. Then $\operatorname{int}M$ maps onto $T$, and $T$ onto $W_\infty$.
\end{proof}

As it must be \cite{FN}, the Julia set is uniformly perfect. A ring domain in $\F(f)$ that separates $\J(f)$ cannot lie in $W_\infty$ or in $T$, since each of these has $\J(f)$ in a single complementary component. So it lies in one of the shells $\{e^a<|x|<e^b\}$ and separates its boundary spheres. Normalise the modulus of a ring domain so that the round shell $\{e^a<|x|<e^b\}$ has modulus $b-a$. By the monotonicity of modulus \cite{Ge62}, a ring domain separating the boundary components of a round shell has modulus at most that of the shell. So the modulus is at most $b-a$, which is at most the length $\beta_1-\alpha_2$ of the largest gap of $\Cs$. Since $\alpha_2,\beta_1\in\Cs$, the shell $\operatorname{int}M$ itself separates $\J(f)$ and has modulus $\beta_1-\alpha_2$, so this bound is sharp.

We turn now to the dimension of the Julia set, and then to the dynamics of $f$ on it.

\begin{proposition}\label{prop:dim}
Let $s\in(0,1)$ be the solution of $d^{-s}+m^{-s}=1$. Then $\dim_H\J(f)=2+s$. In particular $\J(f)$ has Lebesgue measure zero.
\end{proposition}

\begin{proof}
The map $\Theta(t,\xi)=e^t\xi$ is a bi-Lipschitz homeomorphism from $[\alpha,\beta]\times\Sph^2$, with the product of the Euclidean metrics, onto the closed shell $\{e^\alpha\le|x|\le e^\beta\}$. By Theorem~\ref{thm:dyn}, $\Theta(\Cs\times\Sph^2)=\J(f)$, and bi-Lipschitz maps preserve Hausdorff dimension. By the proof of Lemma~\ref{lem:1d}, $\Cs$ is the attractor of an iterated function system of two similarities with ratios $1/d$ and $1/m$ and disjoint images. So it satisfies the open set condition, and its Hausdorff and upper box dimensions both equal $s$ \cite{Hu}, \cite[Ch.~9]{Fa}. The product formulas for dimension \cite[Ch.~7]{Fa} give
\[
2+s=\dim_H\Sph^2+\dim_H\Cs\le\dim_H(\Sph^2\times\Cs)\le\dim_H\Sph^2+\overline{\dim}_B\,\Cs=2+s.\qedhere
\]
\end{proof}

To describe $f$ on $\J(f)$, let $\Sigma_2=\{1,2\}^{\N_0}$, with the shift $\sigma$ and the metric $d_{\Sigma_2}(\omega,\omega')=2^{-\min\{j:\,\omega_j\ne\omega'_j\}}$. Let $\phi=\phi_i$ on $\Cs\cap I_i$, and let $\iota\colon\Cs\to\Sigma_2$ be the itinerary map, with $\iota(t)_j=i$ when $\phi^j(t)\in I_i$. As for any iterated function system of contractions with disjoint images, $\iota$ is a homeomorphism with $\iota\circ\phi=\sigma\circ\iota$ \cite{Hu}, \cite[Ch.~9]{Fa}.

By Remark~\ref{rem:polar}, there are branched covers $\Lambda_1,\Lambda_2$ of $\Sph^2$, of degrees $d^2$ and $m^2$, with $\Lambda_1(\psi(y))=\psi(dy)$ and $\Lambda_2(\psi(y))=\psi(m\bar y)$, such that
\begin{equation}\label{eq:radialPQ}
P(e^t\xi)=e^{\phi_1(t)}\Lambda_1(\xi),\qquad Q(e^t\xi)=e^{\phi_2(t)}\Lambda_2(\xi)\qquad(\xi\in\Sph^2).
\end{equation}
These are maps of Lattès type. Let $\mathbb T=\R^2/2\Z^2$. Since $\psi$ is invariant under the translations in $\Gamma$, it induces a continuous surjection $\psi_{\mathbb T}\colon\mathbb T\to\Sph^2$. Its fibres have at most two points, $[y]$ and $[-y]$, by (Z2). The maps $y\mapsto dy$ and $y\mapsto m\bar y$ induce endomorphisms $\widehat\Lambda_1$ and $\widehat\Lambda_2$ of $\mathbb T$, with $\psi_{\mathbb T}\circ\widehat\Lambda_i=\Lambda_i\circ\psi_{\mathbb T}$. Let $\ell_\psi$ be the image under $\psi$ of normalised Lebesgue measure on $[0,2]^2$, that is, the image under $\psi_{\mathbb T}$ of Haar measure on $\mathbb T$.

\begin{proposition}\label{prop:skew}
Let $\Upsilon\colon\Sigma_2\times\Sph^2\to\Sigma_2\times\Sph^2$ be the skew product $\Upsilon(\omega,\xi)=(\sigma\omega,\Lambda_{\omega_0}(\xi))$, and put $\Xi=\Theta\circ(\iota^{-1}\times\mathrm{id})$.
\begin{enumerate}
\item $\Xi$ is a homeomorphism from $\Sigma_2\times\Sph^2$ onto $\J(f)$, and $f\circ\Xi=\Xi\circ\Upsilon$.
\item Repelling periodic points of $f$ are dense in $\J(f)$. More precisely, the points $x\in\J(f)$ with $f^p(x)=x$ for some $p\ge1$, such that near $x$ the map $f^p$ is topologically conjugate to an expanding similarity near its fixed point, are dense in $\J(f)$.
\item $h_{\rm top}(f|_{\J(f)})=\log(d^2+m^2)=\log\deg f$.
\item Let $\nu$ be the Bernoulli measure on $\Sigma_2$ with weights $\varpi_1=d^2/(d^2+m^2)$ and $\varpi_2=m^2/(d^2+m^2)$. Then $\mu=\Xi_*(\nu\times\ell_\psi)$ is an $f$-invariant mixing Borel probability measure on $\J(f)$ with $h_\mu(f)=\log(d^2+m^2)$. So $\mu$ is a measure of maximal entropy for $f|_{\J(f)}$.
\end{enumerate}
\end{proposition}

\begin{proof}
Write $D=d^2+m^2$, $q_1=d$ and $q_2=m$.

\emph{(i)} Recall that $\Cs\cap I_1\subset[\beta_1,\beta]$, where $f=P$, and $\Cs\cap I_2\subset[\alpha,\alpha_2]$, where $f=Q$. By \eqref{eq:radialPQ}, $f(e^t\xi)=e^{\phi(t)}\Lambda_i(\xi)$ for $t\in\Cs\cap I_i$. Composing with $\iota$ gives (i).

\emph{(ii)} Let $x_0=h(y_0,t_0)\in\J(f)$, with $t_0\in\Cs$, and let $\varepsilon>0$. Choose $y_0'\notin\Z^2$ with $|y_0'-y_0|<\varepsilon$. For $p\ge1$, let $t\in\Cs$ be the point whose itinerary is the periodic word repeating the first $p$ letters $\omega_0\dots\omega_{p-1}$ of $\iota(t_0)$. Then $\phi^p(t)=t$ and $|t-t_0|\le2^{-p}(\beta-\alpha)$, since both lie in the same cylinder of length $p$. Put $F=F_{\omega_{p-1}}\circ\dots\circ F_{\omega_0}$. For this periodic height $t$, we have $F(y,t)=(Ay,t)$, where $A=qR_0$ with $q=\prod_{j<p}q_{\omega_j}\ge2^p$ and $R_0$ either the identity or $y\mapsto\bar y$. The eigenvalues of $A$ lie in $\{q,-q\}$, so $A-I$ is invertible with $\|(A-I)^{-1}\|\le1/(q-1)$. For $a\in\Z^2$ put $y_a=(A-I)^{-1}(2a)$ and $u_a=(y_a,t)$. Then $F(u_a)=u_a+(2a,0)$, which lies in the $\Gamma$-orbit of $u_a$, so $x_a=h(u_a)$ satisfies $f^p(x_a)=x_a$. Here we used that $f\circ h=h\circ F_{\omega_j}$ near the $j$-th point of the orbit, whose height lies in $\Cs\cap I_{\omega_j}$.

Since $2\Z^2$ is $\sqrt2$-dense in $\R^2$, the points $y_a$ are $\sqrt2/(q-1)$-dense. For large $p$, this is less than $\min\{\varepsilon,\dist(y_0',\Z^2)/2\}$. So we can choose $a$ with $|y_a-y_0'|<\min\{\varepsilon,\dist(y_0',\Z^2)/2\}$, and then $y_a\notin\Z^2$. Then $u_a\notin B_h$. Moreover, no point of the orbit of $u_a$ under $F_{\omega_0},F_{\omega_1},\dots$ lies in $B_h$. Otherwise, since $F_1$ and $F_2$ map $B_h$ into itself, $F(u_a)$ would lie in $B_h$, and so would $u_a$, as $B_h$ is $\Gamma$-invariant. So $h$ is injective near each point of this orbit, and near $x_a$
\[
f^p=h\circ(\gamma^{-1}\circ F)\circ h_a^{-1},
\]
where $h_a^{-1}$ is the inverse of $h$ near $u_a$, and $\gamma\in\Gamma$ is the translation by $(2a,0)$. Now $\gamma^{-1}\circ F$ is a similarity of ratio $q>1$ fixing $u_a$, so $x_a$ is repelling in the stated sense. Finally, $|y_a-y_0|<2\varepsilon$ and $|t-t_0|\le2^{-p}(\beta-\alpha)$, and $h$ is Lipschitz on slabs, so $x_a$ is close to $x_0$ when $\varepsilon$ is small and $p$ is large.

\emph{(iii)} Topological entropy is preserved by conjugacy, so by (i) it suffices to compute $h_{\rm top}(\Upsilon)$. Let $\widehat\Upsilon(\omega,y)=(\sigma\omega,\widehat\Lambda_{\omega_0}(y))$ on $\Sigma_2\times\mathbb T$. The map $\Pi=\mathrm{id}\times\psi_{\mathbb T}$ satisfies $\Pi\circ\widehat\Upsilon=\Upsilon\circ\Pi$, and its fibres have at most two points. So $h_{\rm top}(\Upsilon)\le h_{\rm top}(\widehat\Upsilon)$ \cite[Ch.~7]{Wa}, and $h_{\rm top}(\widehat\Upsilon)\le h_{\rm top}(\Upsilon)$ by Bowen's inequality \cite[Theorem~17]{Bo}, since finite sets carry zero entropy. It remains to show that $h_{\rm top}(\widehat\Upsilon)=\log D$.

Give $\mathbb T$ the quotient metric $|\cdot|_{\mathbb T}$ of the Euclidean metric, and $\Sigma_2\times\mathbb T$ the maximum of the two metrics. For a word $w=w_0w_1\cdots$ and $j\ge0$, put $a_j(w)=\prod_{i<j}q_{w_i}$. If $\omega$ begins with $w$, then $\widehat\Lambda_{\omega_{j-1}}\circ\dots\circ \widehat\Lambda_{\omega_0}$ lifts to a linear similarity of ratio $a_j(w)$. Hence, for a lift $v$ of $y'-y$,
\begin{equation}\label{eq:torusexp}
|\widehat\Lambda_{\omega_{j-1}}\circ\dots\circ \widehat\Lambda_{\omega_0}(y')-\widehat\Lambda_{\omega_{j-1}}\circ\dots\circ \widehat\Lambda_{\omega_0}(y)|_{\mathbb T}\le a_j(w)|v|,
\end{equation}
with equality when $a_j(w)|v|<1$, because a vector of length less than $1$ is the shortest representative of its class modulo $2\Z^2$. Write $d_n$ for the Bowen metric $\max_{0\le j<n}$ of the distances between the $j$-th iterates.

\emph{Upper bound.} Let $0<\varepsilon<1$, and choose $k$ with $2^{-k}<\varepsilon$. For each word $w$ of length $n+k$, fix $\omega^w\in\Sigma_2$ beginning with $w$, and a set $Y_w\subset\mathbb T$ that is $(\varepsilon/a_n(w))$-dense in $\mathbb T$, with $\#Y_w\le c\,a_n(w)^2/\varepsilon^2$ for a constant $c$. Let $(\omega,y)\in\Sigma_2\times\mathbb T$, let $w$ be the first $n+k$ letters of $\omega$, and choose $y'\in Y_w$ with $|y-y'|_{\mathbb T}<\varepsilon/a_n(w)$. For $j<n$, the sequences $\sigma^j\omega$ and $\sigma^j\omega^w$ agree in their first $k+1$ letters, and by \eqref{eq:torusexp} the $j$-th torus coordinates are within $a_j(w)\varepsilon/a_n(w)\le\varepsilon$. So $d_n((\omega,y),(\omega^w,y'))\le\varepsilon$. Hence there is an $(n,\varepsilon)$-spanning set with at most
\[
\sum_{|w|=n+k}\frac{c\,a_n(w)^2}{\varepsilon^2}=\frac{c\,2^k}{\varepsilon^2}\sum_{|w'|=n}\ \prod_{i<n}q_{w'_i}^2=\frac{c\,2^k}{\varepsilon^2}D^n
\]
points, and $h_{\rm top}(\widehat\Upsilon)\le\log D$.

\emph{Lower bound.} Let $0<\varepsilon\le1/(2\max\{d,m\})$. For each word $w$ of length $n$, fix $\omega^w$ beginning with $w$, and let $Y'_w\subset\mathbb T$ be a maximal set whose points are at mutual distance at least $\varepsilon/a_{n-1}(w)$. Comparing areas, $\#Y'_w\ge4a_{n-1}(w)^2/(\pi\varepsilon^2)$. Points $(\omega^w,y)$ and $(\omega^{w'},y')$ with $w\ne w'$ are at $d_n$-distance at least $1$, from the first index at which $w$ and $w'$ differ. Now let $y\ne y'$ in $Y'_w$, and let $v$ be a shortest lift of $y'-y$, so $|v|\ge\varepsilon/a_{n-1}(w)$. Let $j\le n-1$ be the least index with $a_j(w)|v|\ge\varepsilon$. If $j=0$, the torus coordinates already differ by $|v|\ge\varepsilon$. If $j\ge1$, then $a_j(w)|v|<\max\{d,m\}\varepsilon\le\frac12$. So by the equality case of \eqref{eq:torusexp}, the $j$-th torus coordinates are at distance $a_j(w)|v|\ge\varepsilon$. Hence the points $(\omega^w,y)$, for $|w|=n$ and $y\in Y'_w$, are at mutual $d_n$-distance at least $\varepsilon$, so they form an $(n,\varepsilon/2)$-separated set with at least
\[
\sum_{|w|=n}\frac{4a_{n-1}(w)^2}{\pi\varepsilon^2}=\frac{8}{\pi\varepsilon^2}D^{n-1}
\]
points, and $h_{\rm top}(\widehat\Upsilon)\ge\log D$.

\emph{(iv)} Let $\widehat\nu=\nu\times\lambda_{\mathbb T}$ on $\Sigma_2\times\mathbb T$, with $\lambda_{\mathbb T}$ the Haar probability measure. Let $\mathcal P$ be the partition of $\Sigma_2\times\mathbb T$ into the $D$ sets $\{\omega_0=i\}\times\pi(R)$, where $R$ runs over the $q_i^2$ half-open squares of side $2/q_i$ tiling $[0,2)^2$, and $\pi\colon\R^2\to\mathbb T$ is the projection. Then $\widehat\Upsilon$ maps each atom $Z$ of $\mathcal P$ bijectively onto $\Sigma_2\times\mathbb T$, and $\widehat\nu(\widehat\Upsilon(\mathcal E))=D\,\widehat\nu(\mathcal E)$ for Borel $\mathcal E\subset Z$. Indeed, on product sets the shift multiplies $\nu$ by $1/\varpi_i$ and $\widehat\Lambda_i$ multiplies $\lambda_{\mathbb T}$ by $q_i^2$, and $q_i^2/\varpi_i=D$. Summing over the atoms, $\widehat\nu(\widehat\Upsilon^{-1}(\mathcal E))=\widehat\nu(\mathcal E)$, so $\widehat\nu$ is invariant.

Let $\mathcal P_n=\bigvee_{j<n}\widehat\Upsilon^{-j}\mathcal P$. By induction, each of its $D^n$ atoms is mapped by $\widehat\Upsilon^n$ bijectively onto $\Sigma_2\times\mathbb T$, and has measure $D^{-n}$. An atom of $\mathcal P_n$ determines the first $n$ letters $w$ of $\omega$, and its torus coordinates lie in the image of a square of side $2/a_n(w)\le2^{1-n}$. So the diameters of the atoms tend to $0$, $\mathcal P$ is a generator, and the Kolmogorov--Sinai theorem gives $h_{\widehat\nu}(\widehat\Upsilon)=\lim_n\frac1n\log D^n=\log D$ \cite[Ch.~4]{Wa}. The same constant-Jacobian property gives $\widehat\nu(Z\cap\widehat\Upsilon^{-n}\mathcal E)=\widehat\nu(Z)\,\widehat\nu(\mathcal E)$ for every atom $Z$ of $\mathcal P_k$, every $n\ge k$ and every Borel set $\mathcal E$. Since finite unions of such atoms approximate every Borel set, $\widehat\nu$ is mixing.

Now $\Pi_*\widehat\nu=\nu\times\ell_\psi$, which is therefore $\Upsilon$-invariant and mixing. By the variational principle \cite[Ch.~8]{Wa} and (iii), $h_{\nu\times\ell_\psi}(\Upsilon)\le\log D$. Conversely, the relativised variational principle of Ledrappier and Walters \cite{LW}, with the integral of the fibre entropies bounded by their supremum, gives
\[
\log D=h_{\widehat\nu}(\widehat\Upsilon)\le h_{\nu\times\ell_\psi}(\Upsilon)+\sup_z h_{\rm top}(\widehat\Upsilon,\Pi^{-1}(z))=h_{\nu\times\ell_\psi}(\Upsilon),
\]
since the fibres of $\Pi$ are finite. Transporting by $\Xi$ gives (iv).
\end{proof}

We use the following standard identity for Möbius maps $\tau$ of $\Sph^n$. Write $\chi$ for the chordal metric, and $|\tau'(x)|_\chi=|\tau'(x)|(1+|x|^2)/(1+|\tau(x)|^2)$ for the chordal derivative.

\begin{lemma}\label{lem:mobius}
For every Möbius map $\tau$ of $\Sph^n$ and all $x,y\in\Sph^n$,
\[
\chi(\tau(x),\tau(y))^2=|\tau'(x)|_\chi\,|\tau'(y)|_\chi\,\chi(x,y)^2 .
\]
\end{lemma}

\begin{proof}
First consider the Euclidean identity $|\tau(x)-\tau(y)|^2=|\tau'(x)|\,|\tau'(y)|\,|x-y|^2$ for finite $x$, $y$, $\tau(x)$, $\tau(y)$.
\begin{itemize}
\item It holds for similarities.
\item It holds for the inversion $x\mapsto x/|x|^2$, since $\big|x/|x|^2-y/|y|^2\big|=|x-y|/(|x||y|)$ and $|\tau'(x)|=|x|^{-2}$.
\item It is preserved under composition, by the chain rule.
\end{itemize}
Every Möbius map is a composition of similarities and this inversion, so the identity holds in general; see also \cite[(3.1.5) and \S3.5]{Be}, where it is written as $|\tau(x)-\tau(y)|=rt^2|x-y|/(|x-a|\,|y-a|)$ for $\tau(\infty)\ne\infty$. Substituting it into the formula $\chi(x,y)=2|x-y|\big((1+|x|^2)(1+|y|^2)\big)^{-1/2}$ \cite[(3.1.3)]{Be} gives the chordal identity, which extends to points at $\infty$ by continuity.
\end{proof}

\begin{remark}\label{rem:consistency}
Three features of the construction are necessary, in the following sense.
\begin{enumerate}
\item \emph{$f$ branches on $\J(f)$.} The branch rays of $P$ and $Q$ cross every sphere $S_t$ with $t\in \Cs$. Some branching on $\J(f)$ is unavoidable. The sphere $S_\beta$ is invariant, and $f|_{S_\beta}=P|_{S_\beta}$ has degree $d^2\ge4$. If $f$ were locally injective near $S_\beta$, then $f|_{S_\beta}$ would be a local homeomorphism of $S^2$ onto itself, by invariance of domain, hence a covering map of degree $d^2$. This is impossible, as $S^2$ is simply connected.
\item \emph{$f$ is not conformal near $\J(f)$.} Indeed, $f$ is of Lattès type on each sphere $S_t$, $t\in \Cs$. A map that is conformal and uniformly expanding near its Julia set could not have a Julia set with non-degenerate components. To see this, suppose that a uqr map $g$ of $\Sph^n$, $n\ge3$, is conformal on a neighbourhood $U$ of $\J(g)$, with $|g'|_\chi\ge\lambda>1$ there, where $|g'|_\chi$ is the chordal derivative defined before Lemma~\ref{lem:mobius}. Shrinking $U$, we may assume it has finitely many components. By Liouville's theorem \cite{Ge62}, \cite{IM2}, $g$ is a Möbius map on each component of $U$. By Lemma~\ref{lem:mobius}, $g$ expands chordal distances by at least $\lambda$ on each component of $U$, globally and not only infinitesimally. A component $C$ of $\J(g)$ lies in a single component of $U$, and so does each $g^k(C)$. Hence $\diam g^k(C)\ge\lambda^k\diam C$, and $C$ is a point. So $\J(g)$ is totally disconnected and $\F(g)$ is connected.
\item \emph{The shell map is not radial.} Suppose an open map $\Phi'$ on $M$ mapped each sphere $S_t$ into a sphere $S_{b(t)}$. Then $b$ is continuous, and openness forces $b$ to be strictly monotone. Indeed, if $b$ had a local extremum at an interior point $t_0$, strict or not, then a small neighbourhood of a point of $S_{t_0}$ would be mapped into one closed side of $S_{b(t_0)}$, while its image contains a point of $S_{b(t_0)}$; so the image would not be open. But the boundary conditions force $b(\alpha_2)=b(\beta_1)=\alpha$. So the fold in the height dynamics, which produces the Cantor set $\Cs$, cannot be realised by maps preserving concentric spheres. This is why $\Phi$ has to be built by a non-radial construction.
\end{enumerate}
\end{remark}

\begin{remark}[Higher dimensions]\label{rem:highdim}
The construction extends to $\Sph^n$ for every $n\ge3$, giving uqr maps of degree $d^{n-1}+m^{n-1}$ whose Julia set is a Cantor set of round $(n-1)$-spheres, with $\dim_H\J(f)=n-1+s$. We outline the changes. The details, which we have checked, are routine but lengthy, and we omit them.
\begin{itemize}
\item Zorich maps and the semiconjugacies of Section~\ref{sec:zorich} exist in all dimensions \cite[\S6.5.4]{IM2}. Suppose the angular map is obtained by reflecting a bi-Lipschitz map of $[0,1]^{n-1}$ onto a hemisphere in the grid faces. Then its automorphy group consists of the maps $y\mapsto Ey+2a$, where $a\in\Z^{n-1}$ and $E$ is a diagonal sign matrix with determinant $1$. The maps $F_1(y,t)=(dy,\phi_1(t))$ and $F_2(y,t)=(m\hat y,\phi_2(t))$, where $\hat y=(-y_1,y_2,\dots,y_{n-1})$, are orientation-preserving and conjugate this group into itself.
\item In the slit construction, the slit is a codimension-one disc, the new branch set is a codimension-two sphere, and $B^n\#B^n\cong S^{n-1}\times[0,1]$.
\item For $n\ne4$, the map $\Psi$ can be obtained from Theorem~\ref{thm:TV} as in dimension three. In place of Kneser's theorem one uses that every orientation-preserving homeomorphism of $S^{n-1}$ is isotopic to the identity, which follows from the stable homeomorphism theorem, see \cite{Ki69}, \cite{Qu82} and the references there. In every dimension, including $n=4$, $\Psi$ can instead be written down explicitly. Place the slit where the inverse branches $\sigma_P$ and $\sigma_Q$ are smooth. Parametrize the double of a small ellipsoid around the slit by oblate spheroidal coordinates, reparametrized radially near the rim so that the map is bi-Lipschitz in the winding-map charts. Then match the collars using diffeomorphisms supported in small balls, each the time-one map of an explicit vector field.
\end{itemize}
\end{remark}

\section{A no-wandering theorem for expanding uqr maps}\label{sec:nwd}

Throughout this section, $f\colon\Sph^n\to\Sph^n$ is uqr, $n\ge2$, and $\chi$ is the chordal metric. A \emph{wandering orbit} is a sequence $(U_k)_{k\ge0}$ of pairwise distinct Fatou components with $f(U_k)=U_{k+1}$. We fix $z\in U_0$ and write $z_k=f^k(z)$. Branching on $\J(f)$ has to be allowed: by Remark~\ref{rem:consistency}(1), a locally injective expansion hypothesis near $\J(f)$ would exclude the map of Theorem~\ref{thm:main}. So expansion is formulated through inclusions of balls.

\begin{lemma}\label{lem:tojulia}
For every wandering orbit, $\operatorname{dist}_\chi(z_k,\J(f))\to0$.
\end{lemma}

\begin{proof}
For $\varepsilon>0$, the set $\{x:\operatorname{dist}_\chi(x,\J(f))\ge\varepsilon\}$ is a compact subset of $\F(f)$. So it is covered by finitely many Fatou components. Since the $U_k$ are pairwise distinct, it contains $z_k$ for only finitely many $k$.
\end{proof}

\begin{theorem}\label{thm:nwd1}
Suppose there are a neighbourhood $V$ of $\J(f)$ and constants $\delta>0$, $C\ge1$, $\lambda>1$ with the following property. Whenever $k\ge1$, $x\in\Sph^n$ with $x,f(x),\dots,f^{k-1}(x)\in V$, and $0<r\le\delta$,
\[
f^k\big(B_\chi(x,r)\big)\supset B_\chi\big(f^k(x),\ \min\{\lambda^kr,\delta\}/C\big).
\]
Then $f$ has no wandering domains, that is, there is no wandering orbit of Fatou components.
\end{theorem}

\begin{proof}
Let $(U_k)$ be a wandering orbit. By Lemma~\ref{lem:tojulia} there is $k_0$ with $z_k\in V$ for all $k\ge k_0$. Choose $0<r\le\delta$ with $B_\chi(z_{k_0},r)\subset U_{k_0}$. For $k\ge k_0$,
\[
U_k=f^{k-k_0}(U_{k_0})\supset B_\chi\big(z_k,\min\{\lambda^{k-k_0}r,\delta\}/C\big),
\]
which is $B_\chi(z_k,\delta/C)$ for all large $k$. Hence $\operatorname{dist}_\chi(z_k,\J(f))\ge\delta/C$ for all large $k$, contradicting Lemma~\ref{lem:tojulia}.
\end{proof}

The constant $C$ is uniform in $k$. This is the natural form of expansion for maps built from semiconjugacies, where the loss comes only from the semiconjugating map. The hypothesis allows branching on $\J(f)$.

The same expansion hypothesis also controls the periodic Fatou components. Write $\varrho(x)=\dist_\chi(x,\J(f))$.

\begin{lemma}\label{lem:fatouorbits}
Suppose $f$ satisfies the hypothesis of Theorem~\ref{thm:nwd1}. There are $\eta>0$ and $c>0$ such that, for every $y\in\F(f)$:
\begin{enumerate}
\item $\varrho(f^k(y))\ge\eta$ for infinitely many $k$;
\item $\varrho(f^k(y))\ge c\min\{\varrho(y),1\}$ for all $k\ge0$.
\end{enumerate}
\end{lemma}

\begin{proof}
Since $\J(f)$ is compact and $V$ is an open neighbourhood of it, we can choose $\eta$ with $0<\eta<\min\{\delta/C,1\}$ and $W:=\{x:\varrho(x)<\eta\}\subset V$. Then $Y:=\Sph^n\setminus W$ is a compact subset of $\F(f)$.

We first record an estimate. Let $y\in\F(f)$ and $0\le j<k$, and suppose $f^i(y)\in V$ for $j\le i<k$. Put $r=\min\{\varrho(f^j(y)),\delta\}$. The ball $B_\chi(f^j(y),r)$ misses $\J(f)$, so it lies in $\F(f)$, and so does its image under $f^{k-j}$, by forward invariance of $\F(f)$. By hypothesis this image contains $B_\chi\big(f^k(y),\min\{\lambda^{k-j}r,\delta\}/C\big)$. Hence
\begin{equation}\label{eq:Dgrowth}
\varrho(f^k(y))\ge\frac{\min\{\lambda^{k-j}r,\delta\}}{C}\ge\frac rC.
\end{equation}

(i) Suppose $f^i(y)\in W$ for all $i\ge j$. Then \eqref{eq:Dgrowth} applies for every $k>j$. Its middle term equals $\delta/C>\eta$ for large $k$, which contradicts $f^k(y)\in W$.

(ii) If $\F(f)=\emptyset$ there is nothing to prove. Otherwise $Y\ne\emptyset$ by (i), and $f(Y)$ is a compact subset of $\F(f)$. Put $\eta'=\min\{\eta,\min_{f(Y)}\varrho\}>0$, and fix $y\in\F(f)$ and $k\ge0$.
\begin{itemize}
\item If $f^k(y)\in Y$, then $\varrho(f^k(y))\ge\eta\ge\eta'$.
\item If $f^i(y)\in W$ for $0\le i\le k$, then \eqref{eq:Dgrowth} with $j=0$ gives $\varrho(f^k(y))\ge\min\{\varrho(y),\delta\}/C$ when $k\ge1$, and this is trivial when $k=0$.
\item Otherwise, let $j<k$ be the largest index with $f^j(y)\in Y$. Then $\varrho(f^{j+1}(y))\ge\eta'$, and $f^i(y)\in W$ for $j<i\le k$. So $\varrho(f^k(y))\ge\eta'$ if $k=j+1$, and \eqref{eq:Dgrowth}, applied from $j+1$, gives $\varrho(f^k(y))\ge\eta'/C$ if $k>j+1$.
\end{itemize}
Since $C\ge1$ and $\eta'\le\eta<\min\{\delta,1\}$, in all cases
\[
\varrho(f^k(y))\ge\frac{\min\{\varrho(y),\eta'\}}{C}\ge\frac{\eta'}{C}\min\{\varrho(y),1\}.\qedhere
\]
\end{proof}

Recall that a fixed point $x_0$ of a uqr map $g$ is \emph{attracting} if $g^k\to x_0$ uniformly on a neighbourhood of $x_0$, and \emph{superattracting} if, in addition, $x_0$ lies in the branch set of $g$. By \cite[Proposition~4.9]{HMM}, every Fatou component $U$ of $g$ with $g(U)\subset U$ is one of the following:
\begin{enumerate}
\item[(a)] the basin of an attracting or superattracting fixed point in $U$;
\item[(b)] a \emph{parabolic basin}: there are a fixed point $x_0\in\partial U$ and $k_j\to\infty$ with $g^{k_j}\to x_0$ locally uniformly on $U$;
\item[(c)] a \emph{rotation domain}: $g|_U$ is a homeomorphism of $U$, and the group it generates has compact closure in the group of homeomorphisms of $U$, with the topology of locally uniform convergence.
\end{enumerate}
See \cite[Definition~4.8]{HMM} for these definitions, and also \cite{HM}.

\begin{theorem}\label{thm:periodic}
Suppose $f$ satisfies the hypothesis of Theorem~\ref{thm:nwd1}. Then $f$ has only finitely many periodic Fatou components, and each of them contains an attracting or superattracting periodic point. Consequently every Fatou component is eventually mapped into the immediate basin of one of finitely many attracting or superattracting cycles.
\end{theorem}

\begin{proof}
We may assume that $\F(f)\ne\emptyset$. Let $\eta$, $c$ and $Y=\{\varrho\ge\eta\}$ be as in Lemma~\ref{lem:fatouorbits} and its proof.

\emph{Finiteness.} The Fatou components form an open cover of the compact set $Y$ by pairwise disjoint sets, so only finitely many of them meet $Y$. By Lemma~\ref{lem:fatouorbits}(i), every cycle of periodic components has a member that meets $Y$. Distinct cycles have no common members. So there are only finitely many cycles, each of finite length.

\emph{Type.} Let $U$ be a periodic component of period $p$, and put $g=f^p$. Since $f$ maps Fatou components onto Fatou components \cite[Ch.~21]{IM2}, $g(U)=U$, and $g$ is uqr with $\F(g)=\F(f)$. By the classification above, $U$ is of type (a), (b) or (c) for $g$. We exclude (b) and (c).

Suppose $U$ is a parabolic basin. Then there are a fixed point $x_0\in\partial U$ and $k_j\to\infty$ with $g^{k_j}\to x_0$ locally uniformly on $U$. For $y\in U$ we get $\varrho(g^{k_j}(y))\le\chi(g^{k_j}(y),x_0)\to0$, since $x_0\in\J(f)$. This contradicts Lemma~\ref{lem:fatouorbits}(ii), as $\varrho(y)>0$.

Suppose $U$ is a rotation domain, and let $\mathcal G$ be the compact closure of the group generated by $g|_U$. Let $U_i=f^i(U)$ for $0\le i<p$, and put
\[
A=\bigcup_{i=0}^{p-1}\big(U\cap f^{-i}(Y\cap U_i)\big).
\]
Since the Fatou components are open and pairwise disjoint, $Y\cap U_i$ is compact. So $f^{-i}(Y\cap U_i)$ is a compact subset of $\F(f)$, and for the same reason its intersection with $U$ is compact. So $A$ is a compact subset of $U$. Since $U$ is locally compact, the evaluation map $\mathcal G\times U\to U$ is continuous, so $\mathcal G A=\{\gamma(a):\gamma\in\mathcal G,\ a\in A\}$ is compact. But $U$ is a nonempty open proper subset of $\Sph^n$, since $\J(f)\ne\emptyset$ \cite[Ch.~21]{IM2}, so it is not compact. Choose $x\in U\setminus\mathcal G A$.

If $g^k(x)\in A$ for some $k\ge0$, then $x=g^{-k}(g^k(x))\in\mathcal G A$. So $g^k(x)\notin A$ for every $k\ge0$. Now let $N=kp+i$ with $0\le i<p$. Then $f^N(x)=f^i(g^k(x))\in U_i$, and if $f^N(x)\in Y$ then $g^k(x)\in A$. So the $f$-orbit of $x$ never meets $Y$, contradicting Lemma~\ref{lem:fatouorbits}(i).

So $U$ is the basin of an attracting or superattracting fixed point of $g$, which is a periodic point of $f$. For the last statement, let $U'$ be any Fatou component. By Theorem~\ref{thm:nwd1}, $f^j(U')=f^k(U')$ for some $j<k$, so $f^j(U')$ is periodic.
\end{proof}

For a rational map $R$ of degree at least $2$, the hypothesis of Theorem~\ref{thm:nwd1} holds if $R$ is hyperbolic. In that case $R$ is expanding for a suitable conformal metric on a neighbourhood $V$ of $\J(R)$ such that $\overline V$ and $R(\overline V)$ miss the closure of the postcritical set \cite[\S19]{Mil}. So the inverse branches of $R^k$ along orbits in $V$ extend to balls of a fixed radius, and the Koebe distortion theorem turns the expansion into the required ball inclusions. Conversely, if the hypothesis holds, then $V$ contains no critical point $c$, since $R(B_\chi(c,r))$ lies in a ball of radius $O(r^2)$ about $R(c)$. By Theorem~\ref{thm:periodic}, every critical point is then attracted to an attracting cycle, so $R$ is hyperbolic. So for rational maps, Theorem~\ref{thm:nwd1} is the classical fact, going back to Fatou, that hyperbolic rational maps have no wandering domains; see, for example, \cite{Mil}.

\begin{corollary}\label{cor:csphere}
The map $f$ of Theorem~\ref{thm:main} satisfies the hypothesis of Theorem~\ref{thm:nwd1}.
\end{corollary}

\begin{proof}
Write $t(x)=\log|x|$ for the height of $x\ne0,\infty$, and also $t(u)$ for the last coordinate of $u\in\R^3$, so that $t(h(u))=t(u)$. Recall that $\Cs\cap I_1\subset[\beta_1,\beta]$ and $\Cs\cap I_2\subset[\alpha,\alpha_2]$.

\emph{Choice of constants.} By Proposition~\ref{prop:shell} and the definition of $f$, we have $f=P$ on $\{t\ge\beta_1-\kappa\}$ and $f=Q$ on $\{t\le\alpha_2+\kappa\}$. Let $S$ be the slab $\{\alpha-1\le t(u)\le\beta+1\}$. As noted in Section~\ref{sec:zorich}, $h$ is $L$-BLD on $S$ for some $L\ge1$. Hence, for $u$ with $\alpha-\tfrac12\le t(u)\le\beta+\tfrac12$ and $0<\rho\le\tfrac12$, we have the Euclidean inclusions
\[
h(B(u,\rho))\subset B(h(u),L\rho)\qquad\text{and}\qquad h(B(u,\rho))\supset B(h(u),\rho/L).
\]
The first holds because $h$ is $L$-Lipschitz on convex subsets of $S$, and $B(u,\rho)\subset S$. The second follows by path lifting \cite[II.3]{Ri}. Let $b\in B(h(u),\rho/L)$. The segment $[h(u),b]$ avoids $0$: the BLD lower bound applied to the radial derivative $e^t$ of $h$ gives $1/L\le e^{\alpha-1}$, so $\rho/L\le e^{\alpha-1}/2<|h(u)|$. Lift this segment under $h$ starting at $u$. As long as the lift stays in $S$, it has length less than $L\cdot\rho/L=\rho\le\tfrac12$, since $h$ is $L$-BLD there. So it stays in $B(u,\rho)\subset S$, the whole segment lifts, and $b\in h(B(u,\rho))$. On the compact region $R=\{e^{\alpha-1}\le|x|\le e^{\beta+1}\}$ the chordal and Euclidean metrics are comparable: $c_0^{-1}|a-b|\le\chi(a,b)\le c_0|a-b|$ for $a,b\in R$, for some $c_0\ge1$.

Now fix $\varepsilon$ with
\[
0<\varepsilon<\min\Big\{\frac\kappa2,\ \frac{\beta_1-\alpha_2}4,\ \frac1{2\max\{d,m\}},\ \frac{e^{\alpha-1}}2\Big\},
\]
so small that the chordal ball $B_\chi(y,\varepsilon)$ lies in $R$ for every $y$ with $\alpha-\tfrac12\le t(y)\le\beta+\tfrac12$. This is possible because these points $y$ form a compact subset of the interior of $R$. Since $\varepsilon\le e^{\alpha-1}/2$, the Euclidean ball $B(y,\varepsilon)$ also lies in $R$ for these $y$. Let $H_i=\{u\in\R^3:\dist(t(u),\Cs\cap I_i)<2\varepsilon\}$ for $i=1,2$. These sets are disjoint, and $f\circ h=h\circ F_i$ on $H_i$. Put
\[
V=\{x:\dist(t(x),\Cs)<\varepsilon\},\qquad \delta=\varepsilon,\qquad \lambda=2,\qquad C=c_0^2L^2.
\]
Then $V$ is an open neighbourhood of $\J(f)$.

\emph{Following an orbit.} Let $x\in V$ with $x,f(x),\dots,f^{k-1}(x)\in V$, and write $x=h(u_0)$. For $0\le j<k$, let $i_j\in\{1,2\}$ be the index with $\dist(t(f^j(x)),\Cs\cap I_{i_j})<\varepsilon$. Put $u_{j+1}=F_{i_j}(u_j)$, so that $h(u_j)=f^j(x)$. Let $\theta_j$ be the ratio of the similarity $F_{i_{j-1}}\circ\dots\circ F_{i_0}$, with $\theta_0=1$. Each $F_{i_j}$ has ratio $d$ or $m$, so $(\theta_j)$ is non-decreasing and $\theta_j\ge2^j$.

\emph{Claim.} If $0<s\le\varepsilon/\theta_k$, then $f^j\big(h(B(u_0,s))\big)=h\big(B(u_j,\theta_js)\big)$ for $0\le j\le k$.

This holds for $j=0$. Suppose it holds for some $j<k$. The ball $B(u_j,\theta_js)$ has radius at most $\theta_ks\le\varepsilon$, and its centre has height within $\varepsilon$ of $\Cs\cap I_{i_j}$, so it lies in $H_{i_j}$. Hence $f\circ h=h\circ F_{i_j}$ on it, and
\[
f^{j+1}\big(h(B(u_0,s))\big)=h\big(F_{i_j}(B(u_j,\theta_js))\big)=h\big(B(u_{j+1},\theta_{j+1}s)\big).
\]
This proves the claim. It is the step that requires $s\le\varepsilon/\theta_k$: the intermediate balls must stay in the region where the same branch $P$ or $Q$ applies, and their radii grow by the actual multipliers $d$ or $m$, not merely by $2$.

\emph{Conclusion.} Let $0<r\le\delta$ and put $s=\min\{r/(c_0L),\,\varepsilon/\theta_k\}$. The heights of $u_0$ and of $u_k$ are within $\max\{d,m\}\varepsilon\le\tfrac12$ of $[\alpha,\beta]$, and $s,\theta_ks\le\varepsilon\le\tfrac12$. By the choice of $\varepsilon$, the balls below lie in $R$, where the metrics are comparable. So
\[
h(B(u_0,s))\subset B(x,Ls)\subset B\big(x,r/c_0\big)\subset B_\chi(x,r),
\]
and, by the claim,
\[
f^k\big(B_\chi(x,r)\big)\supset h\big(B(u_k,\theta_ks)\big)\supset B\big(f^k(x),\theta_ks/L\big)\supset B_\chi\big(f^k(x),\theta_ks/(c_0L)\big).
\]
Finally,
\[
\theta_ks=\min\Big\{\frac{\theta_kr}{c_0L},\,\varepsilon\Big\}\ge\frac{\min\{2^kr,\delta\}}{c_0L}.
\]
Hence $f^k(B_\chi(x,r))\supset B_\chi\big(f^k(x),\min\{\lambda^kr,\delta\}/C\big)$.
\end{proof}

For the map $f$ of Theorem~\ref{thm:main}, the conclusions of Theorems~\ref{thm:nwd1} and~\ref{thm:periodic} are of course already contained in Theorem~\ref{thm:dyn}. The role of Corollary~\ref{cor:csphere} is to show that the class in Theorem~\ref{thm:nwd1} contains maps with infinitely many Fatou components.

\subsection*{Acknowledgements}
Parts of this work, including the exploration of possible constructions and the drafting of some arguments, were carried out with the assistance of Claude, an AI system developed by Anthropic. The authors have checked all arguments and take full responsibility for the content.

\end{document}